\documentclass[11pt]{article}

\usepackage{amssymb,latexsym}
\usepackage{epsfig}
\usepackage{eufrak}
\usepackage{amsmath}
\usepackage{mathrsfs}
\usepackage{color}

\newtheorem{theorem}{Theorem}
\newtheorem{lemma}{Lemma}
\newtheorem{corollary}{Corollary}

\newcommand{\be}{\begin{equation}}
\newcommand{\ee}{\end{equation}}
\newcommand{\bee}{\begin{eqnarray*}}
\newcommand{\eee}{\end{eqnarray*}}
\newcommand{\bel}{\begin{eqnarray}}
\newcommand{\eel}{\end{eqnarray}}
\newcommand{\bec}{\begin{cases}}
\newcommand{\eec}{\end{cases}}
\newcommand{\bem}{\begin{bmatrix}}
\newcommand{\eem}{\end{bmatrix}}

\newcommand{\la}{\label}
\newcommand{\li}{\left}
\newcommand{\ri}{\right}

\newcommand{\ovl}{\overline}
\newcommand{\udl}{\underline}

\newcommand{\ep}{\epsilon}
\newcommand{\vep}{\varepsilon}

\newcommand{\de}{\delta}

\newcommand{\al}{\alpha}

\newcommand{\f}{\frac}

\newcommand{\cd}{\cdots}

\newcommand{\qu}{\quad}
\newcommand{\qqu}{\qquad}

\newcommand{\bb}{\mathbb}

\newcommand{\wh}{\widehat}

\newcommand{\bs}{\boldsymbol}

\newcommand{\tx}{\text}

\newcommand{\bed}{\begin{description}}
\newcommand{\eed}{\end{description}}
\newcommand{\bei}{\begin{itemize}}
\newcommand{\eei}{\end{itemize}}
\newcommand{\ben}{\begin{enumerate}}
\newcommand{\een}{\end{enumerate}}

\newcommand{\beL}{\begin{lemma}}
\newcommand{\eeL}{\end{lemma}}
\newcommand{\beT}{\begin{theorem}}
\newcommand{\eeT}{\end{theorem}}
\newcommand{\sect}{\section}

\newcommand{\bpf}{\begin{pf}}
\newcommand{\epf}{\end{pf}}
\newcommand{\bsk}{\bigskip}

\newcommand{\pfbox}{\hfill\mbox{$\Box$}}

\newenvironment{pf}{\paragraph*{Proof{\rm.}}}{\pfbox\bigskip}

\begin{document}

\title{{\bf A Link between Binomial Parameters and Means of Bounded Random Variables}
\thanks{The author is currently with Department of Electrical Engineering,
Louisiana State University at Baton Rouge, LA 70803, USA, and Department of Electrical Engineering, Southern University and A\&M College, Baton
Rouge, LA 70813, USA; Email: chenxinjia@gmail.com}}

\author{Xinjia Chen}

\date{February, 2008}

\maketitle

\begin{abstract}

In this paper, we establish a fundamental connection between binomial parameters and means of bounded random variables.  Such connection finds
applications in statistical inference of  means of bounded variables.

\end{abstract}

\section{A Fundamental Identity}

Statistical inference of means of bounded random variables is frequently encountered in many areas of sciences and engineering.  Since there
exist rich techniques for the statistical inference of probability of events based on binomial trial models, it is useful to express  means of
bounded random variables in terms of probability of events.  In this regard, we have derived a fundamental identity in the following theorem.

\begin{theorem} \label{thm1}
Let $X$ be a bounded random variable such that  $0 \leq X \leq c$.  Let $Y$ be a bounded random variable such that  $a \leq Y \leq b$.  Let $U$
be a random variable with a uniform distribution over $[0,1]$. Suppose that $U$ is independent with $X$ and $Y$.  Then,  {\small \[
 \bb{E} [X] =
\bb{E} [Y] + \de  \Pr \{ X \geq Y + \de  U \}  +  |c - a - \de | \Pr \{ X > Y + \de + |c- a - \de |  U \}  -  b  \Pr \{ X < b  U < Y \} \] } for
any real number $\de$.

\end{theorem}

\bsk

 In Theorem \ref{thm1}, if we set $Y = c_1$ and $\de = c_2 - c_1$ with constants $c_1$ and $c_2$, then we have

\begin{corollary} \label{thm2}
Let $X$ be a bounded random variable such that $0 \leq X \leq c$.  Let $0 \leq c_1 \leq c_2 \leq c$.  Let $U$ be a random variable with a
uniform distribution over $[0,1]$. Suppose that $X$ and $U$ are independent.  Then,  \[
 \bb{E} [X] =  c_1 + (c_2 -c_1)  \Pr \{ X \geq c_1 + (c_2 -
c_1) U \} +  (c-c_2)  \Pr \{  X > c_2 + (c-c_2) U \}  -  c_1  \Pr \{ X < c_1 U \}. \]

\end{corollary}

In Corollary 1, setting $c_1 = 0$ and $c_2 = c$, we have

\begin{corollary} \label{thm3}
Let $X$ be a bounded random variable such that $0 \leq X \leq c$.  Let $U$ be a random variable with a uniform distribution over $[0,1]$.
Suppose that $X$ and $U$ are independent.  Then,
\[
\bb{E} [X]  =  c \Pr \{ X \geq c U  \}.
\]

\end{corollary}

\sect{Applications}

The identity in Theorem \ref{thm1} can be used to construct confidence intervals for means of bounded random variables.  Specifically, for the
problem of estimating the mean of random variable $X \in [0, c]$, we can generate $n$ random tuples $(X_i, Y_i, U_i), \; i = 1, \cd, n$  such
that

\bed

\item [(i)] $X_i$ are i.i.d. random samples of $X$;

\item [(ii)] $Y_i$ are i.i.d. random samples of $Y$;

\item [(iii)] $U_i$ are i.i.d. random samples uniformly distributed over $[0, 1]$;

\item [(iv)]  $U_i$ is independent of any tuple $(X_i, Y_i)$.

\eed

In terms of these samples, we can define Bernoulli random variables
\[
A_i = \bec 1 & \tx{for} \; X_i \geq Y_i + \de  U_i,\\
0 & \tx{otherwise} \eec
\]
\[
B_i = \bec 1 & \tx{for} \; X_i > Y_i + \de + |c- a - \de |  U_i ,\\
0 & \tx{otherwise}  \eec
\]
\[
C_i = \bec 1 & \tx{for} \; X_i < b  U_i <
Y_i,\\
0 & \tx{otherwise} \eec
\]
for $i = 1, \cd, n$.  Clearly, $A_i,  \; i = 1, \cd, n$ are i.i.d. Bernoulli random variables with parameter
\[
p_a = \Pr \{ X \geq Y + \de  U \};
\]
$B_i,  \; i = 1, \cd, n$ are i.i.d. Bernoulli random variables with parameter
\[
p_b = \Pr \{ X > Y + \de + |c- a - \de |  U \};
\]
and $C_i,  \; i = 1, \cd, n$ are i.i.d. Bernoulli random variables with parameter
\[
p_c = \Pr \{ X < b  U < Y \}.
\]

Now define
\[
K_a = \sum_{i = 1}^n A_i, \qqu K_b = \sum_{i = 1}^n B_i, \qqu K_c = \sum_{i = 1}^n C_i.
\]
Then, $K_a, \; K_b, \; K_c$ are binomial random variables with parameters $p_a, \; p_b$ and $p_c$ respectively.  Applying the identity in
Theorem \ref{thm1}, we have
\[
\bb{E} [X] = \bb{E} [Y] + \de \; p_a   +  |c - a - \de | \; p_b   -  b \; p_c,
\]
which is the desired link between binomial parameters and means of bounded random variables.  Assuming that $\bb{E} [Y]$ is available without
sampling, we can estimate $\bb{E} [X]$ by means estimating $p_a, \; p_b$ and $p_c$. Clearly, the minimum variance unbiased estimates of $p_a, \;
p_b, \; p_c$ are respectively
\[
\f{ K_a } { n }, \qu \f{ K_b } { n }, \qu \f{ K_c } { n }.
\]
Accordingly, an unbiased estimate of $\bb{E} [X]$ is
\[
\bb{E} [Y] + \de \; \f{ K_a } { n } +  |c - a - \de | \f{ K_b } { n } - b \f{ K_c } { n }.
\]

To construct a confidence interval for $\bb{E} [ X ]$ with confidence coefficient $1 - \al$ where $0 < \al < 1$, we can first construct three
confidence intervals employing the classical method of Clopper and Pearson \cite{CP}  such that
\[
\Pr \{  \udl{\bs{p}}_a  \leq p_a \leq \ovl{\bs{p}}_a  \} \geq 1 - \f{\al}{3}, \qu  \Pr \{  \udl{\bs{p}}_b  \leq p_b \leq \ovl{\bs{p}}_b  \} \geq
1 - \f{\al}{3}, \qu \Pr \{  \udl{\bs{p}}_c  \leq p_c \leq \ovl{\bs{p}}_c  \} \geq 1 - \f{\al}{3}.
\]
Then, we can apply Bonferroni's inequality to obtain a confidence interval $[\udl{\bs{\mu}}, \ovl{\bs{\mu}}]$ such that
\[
\udl{\bs{\mu}} = \bb{E} [Y] + \de \udl{\bs{p}}_a +  |c - a - \de | \udl{\bs{p}}_b - b \ovl{\bs{\mu}}, \qqu \ovl{\bs{\mu}} = \bb{E} [Y] + \de
\ovl{\bs{p}}_a +  |c - a - \de | \ovl{\bs{p}}_b - b \udl{\bs{\mu}}
\]
and
\[
\Pr \{ \udl{\bs{\mu}} \leq  \bb{E} [ X ] \leq \ovl{\bs{\mu}} \} \geq 1 - \al.
\]
Here we have assumed $\de > 0, \; b > 0$ in defining $\udl{\bs{\mu}}$ and $\ovl{\bs{\mu}}$.  The construction of confidence interval is similar
for other signs of $\de$ and $b$.

In a similar spirit, we can construct confidence interval for $\bb{E} [ X ]$ by using the identity in Corollary 1 and Bonferroni's inequality.
In that situation, we need to generate $n$ tuples $(X_i, U_i), \; i = 1, \cd, n$ such that \bed

\item [(i)] $X_i$ are i.i.d. random samples of $X$;

\item [(ii)] $U_i$ are i.i.d. random samples uniformly distributed over $[0, 1]$;

\item [(iii)]  $U_i$ is independent of any $X_i$.

\eed

When using the identity in Corollary 2 to construct confidence interval for $\bb{E} [ X ]$, we don't need to apply Bonferroni's inequality.  The
confidence interval can be obtained by scaling the confidence interval of the corresponding binomial parameter.

Applying the identity in Corollary 2, we can also solve the following sample size problem: \bsk

{\it Suppose $X$ is a random variable bounded in $(0, c)$. How large the sample size should be to ensure that the estimate of $\bb{E} [ X ]$ is
close to $\bb{E} [ X ]$ within error margin $\vep$ with a confidence level at least $1 - \al$?}

By using the link we established at above, we can determine the sample size as follows.

\bed

\item [(i)] Determine the minimum sample size $N$ such that the minimum variance unbiased estimate of any Bernoulli parameter is close to its true value
within error margin $\ep = \f{\vep}{c}$ with a confidence level at least $1 - \al$.  An exact method has recently been developed by Chen
\cite{Chen}.

\item [(ii)] Generate $n$ tuples $(X_i, U_i), \; i = 1, \cd, n$ such that $X_i$ are i.i.d. random samples of $X$;
$U_i$ are i.i.d. random samples uniformly distributed over $[0, 1]$; and $U_i$ is independent of any $X_i$.

\item [(iii)] Return
\[
\wh{ \bs{\mu} } = \f{ c \sum_{i = 1}^N A_i  } { N } \qu \tx{with} \qu  A_i = \bec 1 & \tx{for} \; X_i \geq c  U_i,\\
0 & \tx{otherwise} \eec
\]
as the estimate for $\mu = \bb{E} [ X ]$.

\eed

Then, the estimate guarantees
\[
\Pr \{ | \wh{ \bs{\mu} } - \mu | < \vep \} \geq 1 - \al.
\]

Finally, we would like to point out that, by using the link between the binomial parameters and means of bounded random variables, it is
possible to transform the hypothesis testing problem of a bounded variable mean as the problem of testing a binomial proportion.

\sect{Conclusion}

In this paper, an inherent connection between binomial parameters and means of bounded random variables is established. Such connection can be
useful for the estimation and hypothesis testing of the means of bounded random variables.

 \sect{Proof of Theorem 1}

Throughout the proof, we shall use  Riemann-Stieltjes integration.   Let $F_{X, Y}(.)$ be the jointed cumulative distribution function of random
variables $X$ and $Y$. Let $Z = Y + \de$ and $v = a + \de$.   We need some preliminary results.

\beL \la{lem1}
\[
\int_{ y \leq x \leq z } \;\;x \; d F_{X, Y}  =  \bb{E} [ Y] + \de \; \Pr \{ Z \geq X \geq Y + \de \; U \} - \int_{ x > z } \;\;y \; d F_{X, Y}
- \int_{ x < y } \;\;y \; d F_{X, Y}.
\]
\eeL

\bpf

We claim that \be \la{cm}
 \int_{ y \leq x \leq z } \;\;x \; d F_{X, Y}  =  \de \; \Pr \{ Z \geq X \geq Y + \de \; U \} + \int_{ y \leq x \leq
z } \;\;y \; d F_{X, Y}. \ee

 To show the claim, we need to consider three cases: $\de =0, \; \de < 0$ and $\de > 0$.

In the case of $\de =0$, we have \[ Z = Y, \qu \{ Y \leq X \leq Z \} = \{ Y \leq X \leq Y \} = \{ Y = X  \}, \qu \de \; \Pr \{ Z \geq X \geq Y +
\de \; U \}  = 0. \] Hence,
\[ \int_{ y \leq x \leq z
} \;\;x \; d F_{X, Y}  =  \int_{ y = x  } \;\;x \; d F_{X, Y} = \int_{ y = x  } \;\;y \; d F_{X, Y} = \int_{ y \leq x \leq z } \;\;y \; d F_{X,
Y}
\]
and the claim (\ref{cm}) is true for $\de = 0$.

In the case of $\de < 0$, we have \[ Z < Y, \qqu \{ Y \leq X \leq Z \} = \emptyset,
\]
\[
0 \leq \Pr \{ Z \geq X \geq Y + \de \; U \} \leq \Pr \{ Z \geq Y + \de \; U \} = \Pr \{ \de  \geq  \de \; U \} = \Pr \{  U \geq 1 \} = 0.
\]
and
\[
\int_{ y \leq x \leq z } \;\;x \; d F_{X, Y}  = \int_{ y \leq x \leq z } \;\;y \; d F_{X, Y} = 0.
\]
It follows that the claim (\ref{cm}) is true for $\de < 0$.

In the case of $\de > 0$, using the assumption that $U$ is independent with $X$ and $Y$, we have
\begin{eqnarray*}
\de \; \Pr \{ Z \geq X \geq Y + \de \; U \} &  = & \de \; \Pr \{ Z \geq X \geq  Y, \;\;
Z \geq X \geq Y + \de \; U \}\\
& & = \de \; \int_{y \leq x \leq z, \atop
{z \geq x \geq y + \de \; u , \atop {0 \leq u \leq 1} } }  \;du \; d F_{X, Y} \\
& & = \de \; \int_{ y \leq x \leq z} \li [ \int_{u=0}^{  \min \left( \frac{x - y } { \de},
\;\frac{z - y } { \de},\;1  \right) } du \ri ]  d F_{X, Y} \\
& & = \de \; \int_{ y \leq x \leq z} \; \min \left( \frac{x - y } { \de},
\;\frac{z - y } { \de},\;1  \right)     \; d F_{X, Y}\\
& & = \de \; \int_{ y \leq x \leq z} \;
\frac{x - y } { \de} \; d F_{X, Y}\\
& & = \int_{ y \leq x \leq z } \;\;x \; d F_{X, Y} - \int_{ y \leq x \leq z } \;\;y \; d F_{X, Y}.
\end{eqnarray*}
Rearranging this equation shows that the claim (\ref{cm}) is true for $\de > 0$.

Finally, the proof of the lemma is completed by applying the established claim (\ref{cm}) and the observation that
\[
 \int_{ y \leq x \leq z }\;\;y \; d F_{X, Y}
  =  \bb{E} [ Y]  - \int_{ x > z } \;\;y \; d F_{X, Y} - \int_{ x < y } \;\;y \; d F_{X, Y}.
\]

\epf

\beL \la{lem2}
\[
 \int_{ x > z }
\;\;x \; d F_{X, Y} = |c - v |\; \Pr \{ X
> Z + |c- v | \; U \} \; + \; \int_{ x > z }
\;\;z \; d F_{X, Y}.
\]
\eeL

\bpf

To show the lemma,  we need to consider two cases: $c- v \leq 0$ and $c- v > 0$.

In the case of $c- v \leq 0$, we have $X \leq c \leq a + \de \leq Z$.  Hence,
\[
\{ X > Z \} = \emptyset, \qqu  \{ X
> Z + |c- v | \; U \} = \emptyset,
\]
\[
\int_{ x > z } \;\;x \; d F_{X, Y} = 0 = \int_{ x > z } \;\;z \; d F_{X, Y},  \qqu \Pr \{ X > Z + |c- v | \; U \} = 0.
\]
It follows that the lemma is true for $c- v \leq 0$.

In the case of $c- v > 0$,  note that
\begin{eqnarray*}
|c - v |\; \Pr \{ X
> Z + |c- v | \; U \} &  = & |c - v |\; \Pr \{X
> Z, \;\; X
> Z + |c- v | \; U \}\\
& & = |c - v | \; \int_{ x
> z, \atop{ x
> z + |c- v | \; u, \atop 0 \leq u \leq 1 } } \;  du \; d F_{X, Y}\\
& & = |c - v | \; \int_{ x
> z } \li [
\int_{u=0}^{ \min \left( \frac{ x - z} {|c- v|} ,\; 1 \right)
} \; du  \ri ]  d F_{X, Y} \\
& & = |c - v | \; \int_{ x
> z } \;
\min \left( \frac{ x - z} {|c- v|} ,\; 1 \right)\; d F_{X, Y}\\
& & = |c - v | \; \int_{ x
> z } \;
\frac{ x - z} {|c- v|} \; d F_{X, Y}\\
& & = \int_{ x > z } \;\;x \; d F_{X, Y} \; - \; \int_{ x > z } \;\;z \; d F_{X, Y}.
\end{eqnarray*}
Rearranging the last equation shows that the lemma is true for $c- v > 0$.  This competes the proof of the lemma. \epf

\beL \la{lem3}
\[
\int_{ x < y } \;\;x \; d F_{X, Y} = b \; \Pr \{ Y > X  \geq b \; U \}
\]
\eeL

\bpf

To show the lemma,  we need to consider two cases: $b \leq 0$ and $b > 0$.

In the case of $b \leq 0$, we have $X \geq 0 \geq b \geq Y$.  Hence, $\{ Y > X \} = \emptyset$ and
\[
0 \leq \Pr \{ Y > X  \geq b \; U \} \leq \Pr \{ Y > X \} = 0, \qqu \int_{ x < y } \;\;x \; d F_{X, Y}  = 0.
\]
So the lemma is true for $b \leq 0$.   In the case of $b > 0$,  note that
\begin{eqnarray*}
 b \; \Pr \{ Y > X  \geq
b \; U \} & = & b \; \int_{ b u  \leq x < y, \atop 0 \leq u \leq 1 } \; du \; d F_{X, Y}\\
& & = b \; \int_{ x < y } \li [  \int_{u=0}^
{ \min \left( \frac{x} {b} , \; 1 \right) } \;du  \ri ]  d F_{X, Y}\\
& & =  b \; \int_{ x < y } \; \min \left( \frac{x} {b} , \; 1 \right) \;
d F_{X, Y}\\
& & =  b \; \int_{ x < y } \; \frac{x} {b} \;
d F_{X, Y}\\
& & = \int_{ x < y } \;\;x \; d F_{X, Y}
\end{eqnarray*}
for $b > 0$.   This competes the proof of the lemma.

\epf

\beL \la{lem4}
\[
\int_{ x < y } \;\;y \; d F_{X, Y} = b \; \Pr \{ Y > X , \;\; Y > b \; U \}
\]
\eeL

\bpf

To show the lemma,  we need to consider two cases: $b \leq 0$ and $b > 0$.

In the case of $b \leq 0$, we have $X \geq 0 \geq b \geq Y$.  Hence, $\{ Y > X \} = \emptyset$ and
\[
0 \leq \Pr \{ Y > X , \;\; Y > b \; U \} \leq \Pr \{ Y > X \} = 0, \qqu \int_{ x < y } \;\;x \; d F_{X, Y}  = 0.
\]
So the lemma is true for $b \leq 0$.   In the case of $b > 0$, note that
\begin{eqnarray*}
b \; \Pr \{ Y > X , \;\; Y > b \; U \} & = & b \; \int_{ x < y, \atop {y > b \; u, \atop {0 \leq u \leq 1} } }  du \; d F_{X, Y}\\
& & = b \; \int_{ x < y} \li [ \int_{u=0}^{ \min
\left( \frac{y}{b}, \; 1 \right) } du \ri ] d F_{X, Y} \\
& & = b \; \int_{ x < y} \;\min \left( \frac{y}{b}, \; 1 \right)\;
 d F_{X, Y}\\
& & = b \; \int_{ x < y} \;\frac{y}{b} \;
 d F_{X, Y} = \int_{ x < y} \;y \; d F_{X, Y}.
\end{eqnarray*}
This competes the proof of the lemma. \epf

\bsk

Now we are in a position to prove Theorem 1.   Observing that
\begin{eqnarray*}
\bb{E} [X] & = & \int_{ y \leq x \leq z } \;x \; d F_{X, Y} + \int_{ x > z } \;x \; d F_{X, Y} + \int_{ x < y } \;x \; d F_{X, Y}
\end{eqnarray*}
and using Lemmas \ref{lem1}, \ref{lem2}, \ref{lem3}, we have
\begin{eqnarray*}
\bb{E} [X] &  = & \bb{E} [ Y] + \de \; \Pr \{ Z \geq X \geq Y + \de \; U \} + \int_{ x > z }
\;\;(z - y) \; d F_{X, Y} \\
& & + |c - v | \Pr \{ X
> Z + |c- v |  U \} +  b  \Pr \{ Y > X  \geq b  U \} - \int_{ x < y }  y \; d F_{X, Y}.
\end{eqnarray*}
By virtue of
\[
\int_{ x > z } \;\;(z - y) \; d F_{X, Y} = \de \; \Pr \{ X > Z \}
\]
and Lemma \ref{lem4},
 {\small
\bel \bb{E} [X] & = & \bb{E} [ Y] + \de  \Pr \{ Z \geq X \geq Y + \de  U \}  +  \de  \Pr \{ X > Z \}  + |c - v | \Pr \{
X > Z + |c- v |  U \} \nonumber\\
&  &  + \; b  \Pr \{ Y > X  \geq b U \}  -  b  \Pr \{ Y > X, \; Y  > b  U \} \la{eq1}. \eel Noting that
\[
\{ Z \geq X \geq Y + \de \; U \} \cup \{ X > Z \} = \{ X \geq Y + \de \; U  \}
\]
and that
\[
\{ Z \geq X \geq Y + \de \; U \} \cap \{ X > Z \} = \emptyset,
\]
we have \be \la{eq2}
 \de \; \Pr \{ Z \geq X \geq Y + \de \; U \} \; + \; \de \; \Pr \{ X > Z \}  = \de \; \Pr \{ X \geq Y + \de \; U  \}.
\ee Noting that
\[
\{ Y > X  \geq b \; U \} \cap \{ X  < b \; U < Y \} = \emptyset
\]
and that
\[
\{ Y > X  \geq b \; U \} \cup \{ X  < b \; U < Y \} = \{ Y > X, \;\;Y  > b \; U \},
\]
we have \[
 \Pr \{ Y > X  \geq b \; U \} + \Pr \{ X  < b \; U < Y \} = \Pr \{ Y > X, \;\;Y  > b \; U \}.
\]
 Hence, \be \la{eq3}
 b \; \Pr \{ Y > X  \geq b \; U \} \; - \; b \; \Pr \{ Y > X, \;\;Y  > b \; U \}  = - \; b \; \Pr \{ X  < b \; U < Y \}.
\ee Finally, by (\ref{eq1}), (\ref{eq2}) and (\ref{eq3}), we have {\small \[ \bb{E} [X]  =  \bb{E} [ Y] + \de \Pr \{ X \geq Y + \de U  \}  +  |c
- v | \Pr \{ X
> Z + |c- v |  U \}  - b  \Pr \{ X  < b U < Y \}.
\]}
The proof of Theorem 1 is thus completed.

\end{document}